\documentclass[a4paper,12pt]{article}
\usepackage{latexsym}
\usepackage{amsfonts}
\usepackage{xypic}
\newtheorem{theorem}{Theorem}[section]

\title{\sc Survey on preopen sets\thanks{1991 Math.\ Subject
Classification --- Primary: 54-02, 54-06, 54G; Secondary: 54A05,
54A99, 54H05. \protect\newline Key words and phrases --- preopen set,
locally dense set. \protect\newline Research supported partially by the
Japan - Scandinavia Sasakawa Foundation.}}
\author{Julian Dontchev\\Department of Mathematics\\University
of Helsinki\\PL 4, Yliopistonkatu 15\\00014 Helsinki 10\\Finland}
\date{}
\begin{document}
\baselineskip=20pt plus 1pt minus 1pt
\newcommand{\xti}{$(X,\tau,{\cal I})$}
\newcommand{\ij}{${\cal I}|{\cal J}$}
\newcommand{\fxy}{$f \colon (X,\tau) \rightarrow (Y,\sigma)$}
\newcommand{\dpo}{$\delta$-preopen}
\newcommand{\dpc}{$\delta$p-closed}
\newcommand{\dpco}{$\delta$p-compact}
\newcommand{\ptc}{pre-$\theta$-converge}
\newcommand{\pta}{pre-$\theta$-accumulate}
\maketitle
\begin{abstract}
The aim of this survey article is to cover most of recent
research on preopen sets. I try to present majority of the
results on preopen sets that I am aware of.
\end{abstract}

\section{Introduction}\label{s1}

In 1957, J.D.~Weston ({\em J. London Math. Soc.}, {\bf 32}
(1957), 342--354) showed that a continuous bijection $\phi$ on
a complete metric space $X$ to a Hausdorff space $Y$ is open
if and only if it is nearly open, i.e., if and only if $\phi(G)$
is in the interior of its closure for every open $G$. By simple
derivations from this result, open mapping and continuity
theorems in certain spaces were given in 1974 by B.J.~Pettis
({\em Bull. London Math. Soc.}, {\bf 6} (1974), 37--41) and
applied to groups and linear spaces. 

In 1976, T.~Byczkowski and R.~Pol, (On the closed graph and open
mapping theorems, {\em Bull. Acad. Polon. Sci. S'r. Sci. Math.
Astronom. Phys.}, {\bf 24} (9) (1976), 723--726) extended
Weston's theorem to the case of a \v{C}ech-complete space $X$,
showed that `bijection' cannot be replaced by `surjection', and
listed results corresponding to those derived from Weston's.
Although not stated, it follows that a continuous homomorphism
$\phi$ on a topological group $X$ to a Hausdorff group is open
if $X/\text{Ker} \phi$ is \v{C}ech-complete and $\phi$ maps
identity neighborhoods into somewhere dense sets.

The main result of Byczkowski and Pol is that a nearly continuous
mapping from Hausdorff space to \v{C}ech-complete space is
continuous if it has a closed graph and the preimage of every
compact set is compact. The interesting proof implicitly
establishes an open relations theorem slightly more general than
the above extension of Weston's.

In 1974, P.~Mah and S.S.~Naimpally, (Open and uniformly open
relations, {\em Bull. Amer. Math. Soc.}, {\bf 80} (1974),
1157--1159) announced several open mapping theorems for relations
$R \subseteq X \times Y$, where $(X,\delta)$ is an Efremovic
proximity space and $Y$ a topological space. One of their main
results states that, in this setting, $R$ is open if both weakly
open and nearly open. These results are used to generalize some
theorems of J.D.~Weston, J.L.~Kelley (General topology,
Van Nostrand, Toronto, Ont., 1955) and B.J.~Pettis.

Here is the fundamental definition. A subset $A$ of a topological
space $(X,\tau)$ is called {\em preopen} or {\em locally dense}
or {\em nearly open} if $A \subseteq {\rm Int} ({\rm Cl} (A))$.
A set $A$ is called {\em preclosed} if its complement is preopen
or equivalently if ${\rm Cl} ({\rm Int} (A)) \subseteq A$. If $A$
is both preopen and preclosed, then it is called {\em preclopen}.
The term `preopen' was used for the first time by A.S. Mashhour,
M.E. Abd El-Monsef and S.N. El-Deeb ({On pre-continuous and weak
pre-continuous mappings}, {\em Proc. Math. Phys. Soc. Egypt},
{\bf 53} (1982), 47--53), while the concept of a `locally dense'
set was introduced by H.H.~Corson and E.~Michael ({Metrizability
of certain countable unions}, {\em Illinois J. Math.}, {\bf 8}
(1964), 351--360). 

A function \fxy\ is called precontinuous if preimages of open
sets are preopen (A.S. Mashhour, M.E. Abd El-Monsef and S.N.
El-Deeb, same paper as above). Precontinuity was called by
V.~Pt k nearly continuity (Completeness and open mapping theorem,
{\em Bull. Soc. Math. France}, {\bf 86} (1958), 41--74). Nearly
continuity or precontinuity is known also as almost continuity
(T. Husain, {Almost continuous mappings}, {\em Prace Mat.}, {\bf
10} (1966) 1--7). Precontinuity was studied for real-valued
functions on Euclidean spaces by Blumberg back in 1922 (H.
Blumberg, {New properties of all real functions}, {\em Trans.
Amer. Math. Soc.}, {\bf 24} (1922) 113--128).

\vspace{2.5mm}

\section{Characterizations and general properties}\label{s2}

\vspace{2.5mm}

None of the implications in the following diagram is reversible:

$$
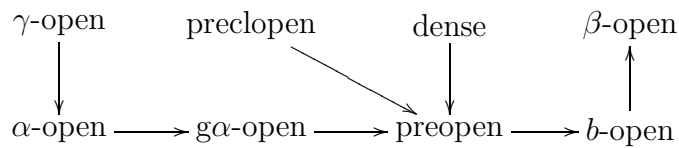
\diagram
\text{$\gamma$-open} \dto & \text{preclopen} \drto & \text{dense}
\dto & \text{$\beta$-open} \\
\text{$\alpha$-open} \rto & \text{g$\alpha$-open} \rto &
\text{preopen} \rto & \text{$b$-open} \uto
\enddiagram
$$
            
\begin{theorem}
For a subset $A$ of a topological space $(X,\tau)$ the following
conditions are equivalent:

{\rm (1)} $A$ is preopen.

{\rm (2)} The semi-closure of $A$ is a regular open set.

{\rm (3)} $A$ is the intersection of an open set and a dense set.

{\rm (4)} $A$ is the intersection of a regular open set and a
dense set.

{\rm (5)} $A$ is a dense subset of some regular open subspace.

{\rm (6)} $A$ is a dense subset of some open subspace.

{\rm (7)} $A$ is a dense subset of some preopen subspace.

{\rm (8)} $A$ is a preneighborhood of each one of its points.

{\rm (9)} ${\rm sCl} (A) = {\rm Int} ({\rm Cl} (A))$.

{\rm (10)} There exists a regular open set $R$ containing $A$
such that ${\rm Cl} (A) = {\rm Cl} (R)$.
\end{theorem}

Here are the most fundamental properties of preopen sets:

$\bullet$ {\em Noiri's Lemma}: If $A$ is semi-open and $B$ is
preopen, then $A \cap B$ is semi-open in $B$ and preopen in $A$.
({T. Noiri}, {Hyperconnectedness and preopen sets}, {\em Rev.
Roumaine Math. Pures Appl.}, {\bf 29} (4) (1984), 329--334.)

$\bullet$ {\em Jankovi\'c and Reilly's Lemma}: Every singleton
is ether preopen or nowhere dense.

$\bullet$ Arbitrary union of preopen sets is preopen.

$\bullet$ Finite intersection of preopen sets need not be
preopen.

$\bullet$ The intersection of a preopen set and an $\alpha$-open
set is a preopen set.

$\bullet$ The intersection $P \cap R$ of a preopen set $P$ and
a regular closed (resp.\ regular open) set $R$ is regular closed
(resp.\ regular open) in the preopen subspace $P$.

$\bullet$ A set is $\alpha$-open if and only if it is semi-open
and preopen.

$\bullet$ A set is clopen if and only if it is closed and
preopen.

$\bullet$ A set is open if and only if it is locally closed and
preopen if and only if it is an $\cal A$-set and preopen if and
only if it is a $\cal B$-set and preopen.

$\bullet$ A set is regular open if and only if it is semi-closed
and preopen.

$\bullet$ If $U$ is a preopen subspace of a space $(X,\tau)$ and
$V$ a preopen subset of $(U,\tau|U)$, then $V$ is preopen in
$(X,\tau)$.

$\bullet$ If $V$ is preopen such that $U \subseteq V \subseteq
{\rm Cl} (U)$, then $U$ is also preopen.

$\bullet$ If $V$ is preopen such that $V \subseteq U \subseteq
X$, then $V$ is also preopen in $(U,\tau|U)$.

$\bullet$ If $A$ is an $\alpha$-open subset of a space
$(X,\tau)$, then a subset $U$ of $A$ is preopen in $(A,\tau|A)$
if and only if $U$ is preopen in $(X,\tau)$.

$\bullet$ If $A$ is a preopen subset of a space $(X,\tau)$, then
for every subset $U$ of $A$ we have $A \cap {\rm sCl} (U) = {\rm
sCl}_{A} (U)$.

$\bullet$ If $A$ is a preopen subset of a space $(X,\tau)$, then
for every subset $U$ of $A$ we have ${\rm Int}_A ({\rm Cl}_A (U))
= A \cap {\rm Int} ({\rm Cl} (U))$.

$\bullet$ If $P$ is preopen and $S$ is semi-open, then $P \cap
{\rm Cl} S = {\rm Cl} (P \cap {\rm Int} (S)) = {\rm Cl} (P \cap
S) = {\rm Cl} (P \cap {\rm Cl} (S)) = {\rm Cl} ({\rm Int} ({\rm
Cl} (P) \cap S))$.

$\bullet$ If $A$ is preopen, then ${\rm Cl} (A) = {\rm
Cl}_{\theta} (A)$.

$\bullet$ $PO(X,\tau) = PO(X,\tau^{\alpha})$.

$\bullet$ If $P$ is preopen in $(X,\tau)$ and $A$ is an
$N$-subset of $(P,\tau|P)$, then $A$ is an $N$-subset of
$(X,\tau)$.

$\bullet$ Let $(X_i)_{i \in I}$ be a family of spaces and
$\emptyset \not= A_i \subseteq X_i$ for each $i \in I$. Then,
$\prod_{i \in I} A_i$ is preopen in $\prod_{i \in I} X_i$ if and
only if $A_i$ is preopen in $X_i$ for each $i \in I$ and $A_i$
is non-dense for only finitely many $i \in I$.

\vspace{2.5mm}

\section{Spaces defined in terms of preopen sets}\label{s3}

\vspace{2.5mm}

The following spaces can be characterized (or defined) in terms
of preopen sets:

A topology $\tau$ on a set $X$ is {\em extremally disconnected}
({M.H. Stone}, {Algebraic characterizations of special Boolean
rings}, {\em Fund. Math.}, {\bf 29} (1937), 223--302) if the
$\tau$-closure of every member of $\tau$ is also in $\tau$.

\begin{theorem}
For a topological space $(X,\tau)$ the following conditions are
equivalent:

{\rm (1)} $X$ is extremally disconnected.

{\rm (2)} Every regular closed subset of $X$ is preopen.

{\rm (3)} Every semi-open subset of $X$ is preopen.

{\rm (4)} The closure of every preopen set is open.

{\rm (5)} The closure of every preopen set is preopen.
\end{theorem}

A nonempty space $X$ is {\em hyperconnected} if every nonempty
open subset of $X$ is dense. A hyperconnected space is called
sometimes {\em irreducible}.

\begin{theorem}
For a topological space $(X,\tau)$ the following conditions are
equivalent:

{\rm (1)} $X$ is hyperconnected.

{\rm (2)} For the semi-closure ${\rm sCl} (A)$ of every nonempty
preopen subset $A$ of $X$ we have ${\rm sCl} (A) = X$, i.e.,
every nonempty preopen subset of $X$ is semi-dense.

{\rm (3)} Every nonempty preopen subset of $X$ is dense.

{\rm (4)} If $S$ is semi-open and $P$ is preopen such that $S
\cap P = \emptyset$, then $S = \emptyset$ or $P = \emptyset$.
\end{theorem}

Recall that a space $(X,\tau)$ is called {\em locally indiscrete}
if every open subset of $X$ is closed.

\begin{theorem}
For a topological space $(X,\tau)$ the following conditions are
equivalent:

{\rm (1)} $X$ is locally indiscrete.

{\rm (2)} Every subset of $X$ is preopen.

{\rm (3)} Every singleton in $X$ is preopen.

{\rm (4)} Every closed subset of $X$ is preopen.
\end{theorem}

\begin{theorem}
For a topological space $(X,\tau)$ the following conditions are
equivalent:

{\rm (1)} $X$ is locally connected.

{\rm (2)} Every component of every open subspace is preopen.
\end{theorem}

Recall that a space $(X,\tau)$ is called {\em submaximal} if
every dense subset of $X$ is open.

\begin{theorem}
For a topological space $(X,\tau)$ the following conditions are
equivalent:

{\rm (1)} $X$ is submaximal.

{\rm (2)} Every preopen set is open.
\end{theorem}

The definition of strongly compact spaces is due to
D.~Jankovi\'{c}, I.~Reilly and M.~Vamanamurthy ({On strongly
compact topological spaces}, {\em Questions Answers Gen.
Topology}, {\bf 6} (1) (1988), 29--40), while the definition of 
strongly Lindel\"of spaces is due to  A.S.~Mashhour, M.E.~Abd.
El-Monsef, I.A.~Hasanein and T.~Noiri ({\em Delta J. Sci.}, {\bf
8} (1984), 30--46). 

\begin{theorem}
For a topological space $(X,\tau)$ the following conditions are
equivalent:

{\rm (1)} $X$ is strongly compact (resp.\ strongly Lindel\"of).

{\rm (2)} Every preopen cover of $X$ has a finite (resp.\
countable) subcover.
\end{theorem}

The definition of countably $P$-compact is due to S.S.~Thakur and
P.~Paik, ({Countably $P$-compact spaces}, {\em Sci. Phys. Sci.},
{\bf 1} (1) (1989), 48--51). 

\begin{theorem}
For a topological space $(X,\tau)$ the following conditions are
equivalent:

{\rm (1)} $X$ is countably $P$-compact.

{\rm (2)} Every countable family of preopen sets which covers $X$
has a finite subcover.
\end{theorem}

The following result is due to M.~Ganster, ({Preopen sets and
resolvable spaces}, {\em Kyungpook Math. J.}, {\bf 27} (2)
(1987), 135--143).

\begin{theorem}
For a connected topological space $(X,\tau)$ the following
conditions are equivalent:

{\rm (1)} $X$ is resolvable.

{\rm (2)} The topology on $X$ having the preopen sets of
$(X,\tau)$ as a subbase is the discrete one.
\end{theorem}

In 1987, V. Popa ({Properties of H-almost continuous
functions}, {\em Bull. Math. Soc. Sci. Math. R.S. Roumanie
(N.S.)}, {\bf 31} (79) (1987), 163--168) introduced the class of
preconnected spaces.

\begin{theorem}
For a topological space $(X,\tau)$ the following conditions are
equivalent:

{\rm (1)} $X$ is preconnected.

{\rm (2)} $X$ cannot be represented as the disjoint union of two
preopen subsets.
\end{theorem}

The definition of strongly irresolvable spaces is due to J.~Foran
and P.~Liebnitz ({A characterization of almost resolvable
spaces}, {\em Rend. Circ. Mat. Palermo}, Serie II, Tomo {\bf XL}
(1991), 136--141). They call a space $X$ {\em strongly
irresolvable} if no nonempty open set is resolvable, where a set
is resolvable if it is resolvable as a subspace. 

\begin{theorem}
For a topological space $(X,\tau)$ the following conditions are
equivalent:

{\rm (1)} $X$ is strongly irresolvable.

{\rm (2)} Every preopen subset is semi-open.

{\rm (3)} Every preopen subset is $\alpha$-open.
\end{theorem}

The concept of a p-closed space is due to Abd El-Aziz Ahmed
Abo-Khadra, {On generalized forms of compactness}, Master's
Thesis, Faculty of Science, Tanta University, Egypt, 1989.

\begin{theorem}
For a topological space $(X,\tau)$ the following conditions are
equivalent:

{\rm (1)} $X$ is p-closed.

{\rm (2)} For every preopen cover $\{ V_{\alpha} \colon \alpha
\in A \}$ of $X$, there exists a finite subset $A_0$ of $A$ such
that $X = \cup \{ {\rm pcl}(V_{\alpha}) \colon \alpha \in A_0
\}$.

{\rm (3)} every maximal filter base \ptc s to some point of $X$,

{\rm (4)} every filter base \pta s at some point of $X$,

{\rm (5)} for every family $\{ V_{\alpha} \colon \alpha \in A \}$
of preclosed subsets such that $\cap \{ V_{\alpha} \colon \alpha
\in A \} = \emptyset$, there exists a finite subset $A_0$ of $A$
such that $\cap \{ {\rm pint} (V_{\alpha}) \colon \alpha \in A_0
\} = \emptyset$.
\end{theorem}

Recently, p-closed spaces were extensively studied by me,
M.~Ganster and T.~Noiri. Here are some typical results of our
paper:

\begin{theorem}
{\rm (1)} Let $(X,\tau)$ be QHC and strongly irresolvable. Then
$(X,\tau)$ is p-closed.

{\rm (2)} Let  $(X,\tau)$ be a p-closed $T_0$ space. Then
$(X,\tau)$ is strongly irresolvable.

{\rm (3)} If a topological space $(X,\tau)$ is p-closed and
$\aleph_0$-extremally disconnected (resp.\ extremally
disconnected), then $(X,\tau)$ is nearly compact (resp.\
$s$-closed).
\end{theorem}

Weak separation axioms in terms of preopen sets were introduced
by A.~Kar and P.~Bhattacharyya, ({\em Bull. Cal. Math. Soc.},
{\bf 82} (1990), 415--422) by replacing the word open with
preopen in the classical definitions. Using their technique one
can of course define {\em pre-$T_{\frac{1}{2}}$-spaces} as the
spaces in which every singleton is either preopen or preclosed.
Unfortunatelly, every spaces is a pre-$T_{\frac{1}{2}}$-space as
shown by H.~Maki, J.~Umehara and T.~Noiri, (Every topological
space is pre-$T\sb {1/2}$, {\em Mem. Fac. Sci. K"chi Univ. Ser.
A Math.}, {\bf 17} (1996), 33--42). 

N.~El-Deeb, I.A.~Hasanein, A.S. Mashhour and T.~Noiri (On
$p$-regular spaces, {\em Bull. Math. Soc. Sci. Math. R. S.
Roumanie (N.S.)}, {\bf 27} (75) (1983), no. 4, 311--315)
introduced the class of $p$-regular spaces. A topological space
$X$ is said to be $p$-regular if for each closed subset $F$ and
each point $x \in X \setminus F$ there exist disjoint preopen
sets $U$ and $V$ such that $F \subseteq U$ and $x \in V$. Almost
$p$-regular and $p$-completely regular spaces were defined by
S.R.~Malghan and G.B.~Navalagi (Almost $p$-regular,
$p$-completely regular and almost $p$-completely regular spaces,
{\em Bull. Math. Soc. Sci. Math. R. S. Roumanie (N.S.)}, {\bf 34}
(82) (1990), 317--326).

\vspace{2.5mm}

\section{Functions defined in terms of preopen sets}\label{s4}

Let $V$ (resp.\ $U$) be a preopen subset of $(Y,\sigma)$ (resp.\
$(X,\tau)$). A function $f \colon (X,\tau) \rightarrow
(Y,\sigma)$ is called:

$\bullet$ {\em Precontinuous} if and only if the preimage of
every open subset of $Y$ is preopen in $X$.

$\bullet$ {\em Strongly $M$-precontinuous} (= SMPC) if and only
if the preimage of each preopen subset of $Y$ is open in $X$.

$\bullet$ {\em Preirresolute} if and only if the preimage of
every preopen subset of $Y$ is preopen in $X$.

$\bullet$ {\em Faintly precontinuous} if and only if the preimage
of every $\theta$-open subset of $Y$ is preopen in $X$.

$\bullet$ {\em Strongly faintly precontinuous} if and only if the
preimage of every preopen subset of $Y$ is $\theta$-open in $X$.

$\bullet$ {\em Weakly continuous} if and only if $f^{-1} (V)
\subseteq {\rm Int} f^{-1} ({\rm Cl} (V))$.

$\bullet$ {\em Weakly continuous} if and only if ${\rm Cl}
(f^{-1} (V)) \subseteq f^{-1} ({\rm Cl} (V))$.

$\bullet$ {\em Almost continuous} if and only if $f^{-1} (V)
\subseteq {\rm Int} f^{-1} ({\rm Int} ({\rm Cl} (V)))$.

$\bullet$ {\em Almost open} if and only if $f(U) \subseteq {\rm
Int} (f({\rm Cl} (U)))$.

$\bullet$ {\em Almost weakly open} if and only if $f(U) \subseteq
{\rm Int} ({\rm Cl} (f({\rm Cl} (U))))$.

$\bullet$ {\em $\alpha$-continuous} if and only if $f({\rm Cl}
(U)) \subseteq {\rm Cl} (f(U))$.

$\bullet$ {\em Preclosed} (resp.\ {\em preopen}) if and only if
the image of every closed (resp.\ open) set is preclosed (resp.\
preopen).

\vspace{2.5mm}

\section{Generalizations of preopen sets}\label{s5}

\vspace{2.5mm}

A subset $A$ of a topological space $(X,\tau)$ is called:

$\bullet$ {\em $\cal I$-open} if $A \subseteq {\rm Int}(A^{*})$
(D.~Jankovi\'{c} and T.R.~Hamlett, {Compatible Extensions of
Ideals}, {\em Boll. Un. Mat. It.}, {\bf 7} (1992), 453-465). Note
here that in the case of the minimal ideal, i.e., when ${\cal I}
= \{ \emptyset \}$, the concepts of $\cal I$-open and preopen
sets coincide.

$\bullet$ {\em $\gamma$-open} if $A \cap S$ is preopen, whenever
$S$ is preopen (D.~Andrijevi\'c and M.~Ganster, On PO-equivalent
topologies, IV International Meeting on Topology in Italy
(Sorrento, 1988), {\em Rend. Circ. Mat. Palermo (2) Suppl.}, {\bf
24} (1990), 251--256). Note that two topologies on a set $X$ are
PO-equivalent if their classes of preopen sets coincide.

$\bullet$ {\em preclopen} if $A$ is both preopen and preclosed
(A.~Kar and P.~Bhattacharyya, ({\em Bull. Cal. Math. Soc.},
{\bf 82} (1990), 415--422).

$\bullet$ {\em generalized p-closed} (= gp-closed) if ${\rm pCl}
(A) \subseteq U$ whenever $A \subseteq U$ and $U$ is open (Noiri,
T., Maki, H. and Umehara, J., Generalized preclosed functions,
Mem. Fac. Sci. Kochi Univ. Ser. A Math. {\bf 19} (1998), 13--20).

$\bullet$ {\em regular generalized p-closed} (= rgp-closed) if
${\rm pCl} (A) \subseteq U$ whenever $A \subseteq U$ and $U$ is
regular open (Noiri, Almost p-regular spaces and some functions,
{\em Acta Math. Hungar.}, {\bf 79} (3) (1998), 207--216).

$\bullet$ {\em generalized $\alpha$-closed set} (= {\em
g$\alpha$-closed}) if ${\rm Cl}_{\alpha}(A) \subseteq U$ whenever
$A \subseteq U$ and $U$ is $\alpha$-open ({H. Maki, R. Devi and
K. Balachandran}, {Generalized $\alpha$-closed sets in topology},
{\em Bull. Fukuoka Univ. Ed. Part {\em III}}, {\bf 42} (1993),
13--21).

$\bullet$ {\em $b$-open} or {\em sp-open} if $A \subseteq {\rm
Cl} ({\rm Int} (A)) \cup {\rm Int} ({\rm Cl} (A))$ ({D.
Andrijevi\'{c}}, {On $b$-open sets}, {\em Mat. Vesnik}, {\bf 48}
(1996), 59--64 and {J. Dontchev and M. Przemski}, {On the various
decompositions of continuous and some weakly continuous
functions}, {\em Acta Math. Hungar.}, {\bf 71} (1-2) (1996),
109--120.)

\vspace{2.5mm}

\section{Recent progress in the theory of preopen sets}\label{s6}

\vspace{2.5mm}

The following recent papers (a list obtained from MR), more or
less deal with (to a certain extent) with preopen sets:

$\bullet$ R.H.~Atia, S.N.~El-Deeb and I.A.~Hasanein, A note on
strong compactness and $S$-closedness, {\em Mat. Vesnik} {\bf 6}
(19) (34) (1982), no. 1, 23--28.

$\bullet$ A.S.~Mashhour, I.A.~Hasanein and S.N.~El-Deeb, A note
on semicontinuity and precontinuity, {\em Indian J. Pure Appl.
Math.}, {\bf 13} (10) (1982), 1119--1123. In this note, the
author proves some elementary properties of semicontinuous and
precontinuous functions. 

$\bullet$ A.S.~Mashhour, F.H.~Khedr, I.A.~Hasanein, and
A.A.~Allam, $S$-closedness in bitopological spaces, {\em Ann.
Soc. Sci. Bruxelles S'r. I}, {\bf 96} (2) (1982), 69--76. The
authors show the relationships of pairwise preopen and
$\alpha$-sets to S-closedness. 

$\bullet$ N.~El-Deeb, I.A.~Hasanein, A.S. Mashhour and T.~Noiri,
On $p$-regular spaces, {\em Bull. Math. Soc. Sci. Math. R. S.
Roumanie (N.S.)}, {\bf 27} (75) (1983), no. 4, 311--315. A
topological space $X$ is said to be $p$-regular if for each
closed subset $F$ and each point $x\in X \setminus F$ there exist
disjoint preopen sets $U$ and $V$ such that $F \subseteq U$ and
$x\in V$. The main results of the paper are the following: If a
subspace $S$ of a $p$-regular space is $\alpha$-open, then it is
$p$-regular. An arbitrary product of $p$-regular spaces is
$p$-regular. The image of a regular space by a continuous
preclosed function with compact point inverses is $p$-regular. 

$\bullet$ A.S.~Mashhour, A.A.~Allam, F.S.~Mahmoud and F.H.~Khedr,
On supratopological spaces, {\em Indian J. Pure Appl. Math.},
{\bf 14} (4) (1983), 502--510. Results are established for
supratopological spaces which generalize known results about
semi-open sets, preopen sets, and $\alpha$-open sets. 

$\bullet$ A.S.~Mashhour, M.E.~Abd El-Monsef and I.A.~Hasanein,
On pretopological spaces, {\em Bull. Math. Soc. Sci. Math. R. S.
Roumanie (N.S.)}, {\bf 28} (76) (1984), no. 1, 39--45. A
bijective map $f$: $X\to Y$ is called a prehomeomorphism whenever
$A \subseteq X$ is preopen if and only if $f(A) \subseteq Y$ is
preopen. The paper contains some easy remarks on preopen sets,
prehomeomorphisms and related concepts. In particular, $X$ is
paracompact if and only if every open cover has a locally finite
preopen refinement. 

$\bullet$ T.~Noiri, Hyperconnectedness and preopen sets, {\em
Rev. Roumaine Math. Pures Appl.}, {\bf 29} (4) (1984), 329--334.
This paper gives a large number of new results and improvements
of older results concerning hyperconnected spaces and sets. A
sample of the results: If $A$ is preopen in $X$ and $B$ is semi-
open in $X$, and either $A$ or $B$ is hyperconnected in $X$, then
$A \cap B$ is hyperconnected in $X$. If $X$ is locally
hyperconnected (i.e. every $x\in X$ has a base of open
hyperconnected neighbourhoods) and if $X\sb 0$ is a preopen or
semi-open subset of $X$, then $X\sb 0$ is locally hyperconnected.
If a function $f:X\to Y$ is semicontinuous, then $f(A)$ is
hyperconnected for every preopen hyperconnected subset $A$ of
$X$. 

$\bullet$ B.~Ahmad, and T.~Noiri, The inverse images of
hyperconnected sets, {\em Mat. Vesnik}, {\bf 37} (2) (1985),
177--181. The following results are proved: (1) If $f$ is an
$\alpha$-continuous function and $A$ is a semi-open
hyperconnected set in $X$, then $f(A)$ is hyperconnected. (2) If
$f$ is a semiclosed preserving surjection with preopen
hyperconnected point inverses and $B$ is a preopen
hyperconnected subset of $Y$, then $f\sp {-1}(B)$ is preopen
hyperconnected. 

$\bullet$ S.F.~Tadros, Some special systems of subsets of
topological spaces and separation axioms defined by them, {\em
Soobshch. Akad. Nauk Gruzin. SSR}, {\bf 118} (1) (1985), 53--56.
The first part of the paper contains characterizations of the
families of $\alpha$-open, preopen and $\beta$-open subsets from
an arbitrary topological space. Further the author establishes
some results about the separation properties induced by the
families mentioned, linking them to the classical separation
axioms (for example, characterizations of Hausdorff, regular and
normal spaces in terms of $\alpha$-sets are given).

$\bullet$ V.~Singh, G.~Chae and D.N.~Misra, (Functions with
closed graph and some other related properties, {\em Univ. Ulsan
Rep.}, {\bf 17} (1) (1986), 127--131) studied relationships
between functions (with closed graphs, which are continuous or
nearly continuous, etc.) and separation axioms satisfied by the
range space. A sample theorem of their paper is the following:
If $f\:X\to Y$ is a nearly open (= preopen) surjection with a
closed graph, then $Y$ is Hausdorff. 

$\bullet$ D.~Andrijevi\'c, A note on preopen sets, Third National
Conference on Topology (Italian) (Trieste, 1986), {\em Rend.
Circ. Mat. Palermo (2) Suppl.}, {\bf 18} (1988), 195--201. In
general, the intersection of two preopen sets need not be
preopen, although the union of an arbitrary collection of preopen
sets is preopen. This paper gives a set of necessary
and sufficient conditions for the intersection of two preopen
sets to be preopen.

$\bullet$ M.~Ganster, Preopen sets and resolvable spaces, {\em
Kyungpook Math. J.}, {\bf 27} (2) (1987), 135--143. This paper
provides solutions to the following questions of Katetov. (1)
Obtain conditions under which every dense-in-itself set is
preopen. (2) Find conditions under which the intersection of any
two preopen sets is preopen. (3) Let ${\cal T}\sp *$ denote the
topology on $X$ having the collection of all preopen sets in
$(X,\cal T)$ as a subbase. Find conditions under which $\cal T\sp
*$ is discrete. 

Resolvability is the crucial notion in answering these questions.
A topological space $(X,\cal T)$ is defined to be resolvable if
there is a subset $D$ of $X$ such that $D$ and $X \setminus D$
are both dense in $(X,\cal T)$; otherwise $(X,\cal T)$; is called
irresolvable. 

$\bullet$ D.~Andrijevi\'c, On the topology generated by preopen
sets, {\em Mat. Vesnik}, {\bf 39} (4) (1987), 367--376. The
topology mentioned in the title has recently been introduced in
a joint paper by the Andrijevi\'c and Ganster (same journal {\bf
39} (2) (1987), 115--119). The paper is a continuation of the
study of this new topology with an emphasis on the properties of
its closure operator. 

$\bullet$ A.A.~Allam, A.M.~Zahran and I.A.~Hasanein, On almost
continuous, $S$-continuous and set connected mappings, {\em
Indian J. Pure Appl. Math.}, {\bf 18} (11) (1987), 991--996. In
this paper, the authors strengthen some of the results of
P.E.~Long and D.A.~Carnahan ({\em Proc. Amer. Math. Soc.}, {\bf
38} (1973), 413--418) by using preopenness and openness for both
mappings and sets. For example, the theorem that if
$f\:X\rightarrow Y$ is an open almost continuous mapping then,
for each open $V \subseteq Y$, ${\rm Cl}(f\sp {-1}(V)) \subseteq
f\sp {-1}({\rm Cl}(V))$ is improved by letting $f$ be a preopen
mapping. The same theorem is also established for a preopen set
$V$ with the openness of $f$ dropped. The authors also improve
some of the results of T.~Noiri ({\em J. Korean Math. Soc.}, {\bf
16} (2) (1979/80), 161--166) and ({\em Kyungpook Math. J.}, {\bf
16} (2) (1976), 243--246) dealing with weakly continuous and
almost continuous mappings. 

$\bullet$ V.~Popa, Characterizations of $H$-almost continuous
functions, {\em Glas. Mat. Ser. III}, {\bf 22} (42) (1987), no.
1, 157--161. The author characterizes almost continuous functions
using preopen and preclosed sets. 

$\bullet$ S.S.~Thakur, and P.~Paik, Locally $P$-connected spaces,
{\em J. Indian Acad. Math.}, {\bf 9} (1) (1987), 34--36. A space
$X$ is said to be $P$-connected if it is not the union of two
sets $A$, $B$ with $A \cap {\rm pCl}(B) = \emptyset =
{\rm pCl} (A)\cap B$. A space $X$ is said to be locally
$P$-connected if for every point $x\in X$ and each open
neighbourhood $U$ of $x$ there exists an open $P$-connected set
$G$ such that $x\in G \subseteq U$. The paper contains basic
properties of locally $P$-connected spaces. 

$\bullet$ V.~Popa, Properties of $H$-almost continuous functions,
{\em Bull. Math. Soc. Sci. Math. R. S. Roumanie (N.S.)}, {\bf 31}
(79) (1987), no. 2, 163--168. The basic properties of almost
continuous maps in terms of preopen and preclosed sets are given.

$\bullet$ A.S.~Mashhour, A.A.~Allam, I.A.~Hasanein and
K.M.~Abd-El-Hakeim, On strongly compactness, {\em Bull. Calcutta
Math. Soc.}, {\bf 79} (4) (1987), 243--248. 

$\bullet$ T.~Noiri, Characterizations of extremally disconnected
spaces, {\em Indian J. Pure Appl. Math.}, {\bf 19} (4) (1988),
325--329. In this paper, the author presents several
characterizations of extremally disconnected spaces in terms of
preopen, semi-open and semi-preopen sets. One of these
characterizations says that a space $X$ is extremally
disconnected if and only if the closure of every preopen subset
of $X$ is preopen. Another one asserts that $X$ is extremally
disconnected if and only if for every semi-open set $A \subseteq
X$ and every semi-preopen set $B \subseteq X$, the set $A\cap B$
is semi-open. 

$\bullet$ M.~Ganster, F.~Gressl and I.L.~Reilly, On a
decomposition of continuity, General topology and applications
(Staten Island, NY, 1989), 67--72, Lecture Notes in Pure and
Appl. Math., {\bf 134}, Dekker, New York, 1991. The main result
of this paper is that $f$ is continuous if and only if it is
precontinuous and weakly $\cal B$-continuous. Applications to
topological groups and topological vector spaces follow. 

$\bullet$ S.S.~Thakur and P.~Paik, Countably $P$-compact spaces,
{\em Sci. Phys. Sci.}, {\bf 1} (1) (1989), 48--51. A topological
space $X$ is defined to be countably $P$-compact if every
countable family of preopen sets which covers $X$ has a finite
subcover. This property implies countable compactness but is not
equivalent to it. Characterizations of the notion (similar to
those given for countable compactness) are given in terms of
countable families of preclosed sets having the countable
intersection property and  of precluster points of sequences in
$X$. Countable $P$-compactness is shown to be hereditary for
preopen and preclosed subsets and to be preserved by
precontinuous, $\alpha$-open onto mappings. 

$\bullet$ M.~Ganster, A note on strongly Lindel\"of spaces, {\em
Soochow J. Math.}, {\bf 15} (1) (1989), 99--104. The author
defines $X$ to be $d$-Lindel\"of if each cover of $X$ by dense
subsets has a countable subcover, and proves that a space is
strongly Lindel\"of if and only if it is Lindel\"of and
$d$-Lindel\"of. Maximal strongly Lindel\"of spaces are
investigated. 

$\bullet$ S.S.~Thakur and P.~Paik, On semi $T\sb 0,$ semi $T\sb
1$ and semi $T\sb 2$ spaces, {\em Mathematica (Cluj)}, {\bf 31}
(54) (1989), 91--94. The authors study in detail the semi-$T\sb
0$, semi-$T\sb 1$ and semi-$T\sb 2$ spaces, especially the
relation between these spaces and preopen semicontinuous
functions. 

$\bullet$ I.~Kupka, On some classes of sets related to
generalized continuity, {\em Acta Math. Univ. Comenian.}, {\bf
56/57} (1989), 55--61 (1990). Here the author considers
$\alpha$-sets, semi-open sets and preopen sets. Suitable forms
of generalized homeomorphisms related to these classes of sets
are considered. The main theorem is as follows: If $X,Y$ are
$T\sb 1$ topological spaces, then the classes of
prehomeomorphisms, of
semihomeomorphisms and of $\alpha$-homeomorphisms from $X$ onto
$Y$ coincide.

$\bullet$ I.L.~Reilly and M.K.~Vamanamurthy, On some questions
concerning
preopen sets, {\em Kyungpook Math. J.}, {\bf 30} (1) (1990),
87--93. In the paper of Mashhour, Abd El-Monsef and El-Deeb,
several questions are attributed to Katetov. First, find
necessary and sufficient conditions under which every pre-open
set is open. Second, find conditions under which every
dense-in-itself set is pre-open. Third, find conditions under
which the intersection of any two pre-open sets is pre-open.
Fourth, find conditions under which the topology generated by the
pre-open sets is discrete. The authors solve the first question
by showing that every pre-open set is open if and only if
every dense set is open (a property considered by Bourbaki). They
also provide a partial solution in terms of door spaces. They
show, further, that every subset is pre-open if and only if every
open set is closed (i.e. if there is a base which is a
partition). The authors partly solve the fourth question by
showing that the topology generated by the pre-open sets of any
space which can be partitioned into two disjoint dense sets is
discrete. 

$\bullet$ V.~Popa, Some properties of $H$-almost continuous
multifunctions, {\em Problemy Mat.}, {\bf 10} (1990), 9--26. The
paper presents various characterizations of almost continuous
multivalued maps in the sense of Husain. The results obtained are
formulated in terms of preopen and preclosed sets. 

$\bullet$ M.~Jeli\'c, A decomposition of pairwise continuity,
{\em J. Inst. Math. Comput. Sci. Math. Ser.}, {\bf 3} (1) (1990),
25--29. She showed that a function $f\:(X,\tau\sb 1,\tau\sb 2)\to
(Y,\sigma\sb 1,\sigma\sb 2)$ is pairwise continuous if and only
if it is pairwise precontinuous and pairwise LC-continuous.    

$\bullet$ M.K.~Singal and N.~Prakash, Fuzzy preopen sets and
fuzzy preseparation axioms, {\em Fuzzy Sets and Systems}, {\bf
44} (2) (1991), 273--281. In this paper, the concept of preopen
sets is generalised to the fuzzy setting. Further fuzzy
separation axioms are introduced and investigated with the help
of fuzzy preopen sets.

$\bullet$ Abdulla S.~Bin Shahna, On fuzzy strong semicontinuity
and fuzzy precontinuity, {\em Fuzzy Sets and Systems}, {\bf 44}
(2) (1991), 303--308. Using the notion of fuzzy sets, Azad, in
1981, introduced and studied the concept of semiopen (semiclosed)
sets, semicontinuous mappings, almost continuous mappings, and
weakly continuous mappings in the fuzzy setting. In the same
spirit, defining fuzzy $\alpha$-open ($\alpha$-closed) sets, and
fuzzy preopen (preclosed) sets, the authors introduce and make
a preliminary study of fuzzy strongly semicontinuous, and fuzzy
precontinuous mappings in this paper.

$\bullet$ D.~Andrijevi\'c, On SPO-equivalent topologies, V
International Meeting on Topology in Italy (Italian) (Lecce,
1990/Otranto, 1990), {\em Rend. Circ. Mat. Palermo (2) Suppl.},
{\bf 29} (1992), 317--328. Two topologies on a set $X$ are
SPO-equivalent if their classes of semi-preopen sets $A$ 
coincide. If $\tau$ is a topology on $X$, the topology $\tau\sb
j=\{A \subseteq X \colon \ A\cap B$ is $\tau$-preopen whenever
$B$ is $\tau$-preopen$\}$ is used to identify SPO-equivalent
topologies on $X$; the class of all topologies on $X$ which are
SPO-equivalent to $\tau$ has a largest topology $\tau\sb R$. 

$\bullet$ A.~Kar and P.~Bhattacharyya, Bitopological preopen
sets, precontinuity and preopen mappings, {\em Indian J. Math.},
{\bf 34} (3) (1992), 295--309. The concepts of preopen sets,
precontinuity and preopen mappings in a bitopological space are
introduced in this paper. The conditions under which the various
properties enjoyed by the above concepts in a single topological
space can be generalized to a bitopological space are
investigated.

$\bullet$ S.S.~Thakur and U.D.~Tapi, Set $P$-connected mappings,
{\em 
Mathematica (Cluj)}, {\bf 34} (57) (2) (1992), 183--185. 
A space $X$ is said to be $P$-connected between $A$ and $B$,
where $A,B \subseteq X$, if there exists no preclosed
preopen $F \subseteq X$ such that $A \subseteq F$ and $F\cap
B=\emptyset$. A mapping $f\colon X\to Y$ is said to be set
$P$-connected if for all $A,B \subseteq X$, $f(X)$ is
$P$-connected between $f(A)$ and $f(B)$ when $X$ is $P$-connected
between $A$ and $B$. The authors establish several properties of
set $P$-connected mappings. 

$\bullet$ D.S.~Jankovi\'c, and Ch.~Konstadilaki-Savvopoulou, On
$\alpha$-continuous functions, {\em Math. Bohem.}, {\bf 117} (3)
(1992), 259--270. 

$\bullet$ Asit Kumar Sen and P.~Bhattacharyya, On preclosed
mapping, {\em Bull. Calcutta Math. Soc.}, {\bf 85} (5) (1993),
409--412. 

$\bullet$ S.~Raychaudhuri and M.N.~Mukherjee, On $\delta$-almost
continuity and $\delta$-preopen sets, {\em Bull. Inst. Math.
Acad. Sinica}, {\bf 21} (4) (1993), 357--366. 

$\bullet$ Abdulla S.~Bin Shahna, Mappings in fuzzy topological
spaces, {\em Fuzzy Sets and Systems}, {\bf 61} (2) (1994),
209--213. In this paper, some characterizations of fuzzy strongly
semicontinuous and fuzzy precontinuous mappings are given. The
definitions of $\alpha$-open and preopen mappings are extended
to fuzzy sets.

$\bullet$ T.~Aho and T.~Nieminen, Spaces in which preopen subsets
are semiopen, {\em Ricerche Mat.}, {\bf 43} (1) (1994), 45--59.
A topological space is called a PS-space if every preopen subset
is semiopen; it is called irresolvable if two dense subsets
cannot be disjoint. The authors study several equivalent    
conditions concerning these notions. 

$\bullet$ R.~Mahmoud and D.~Rose, A note on submaximal spaces and
SMPC functions, {\em Demonstratio Math.}, {\bf 28} (3) (1995),
567--573. In this paper, the authors give several
characterizations of submaximal spaces involving SMPC functions. 

$\bullet$ Manindra Chandra Pal and P.~Bhattacharyya, Faint
precontinuous functions, {\em Soochow J. Math.}, {\bf 21} (3)
(1995), 273--289. A function $f\colon X\to Y$, where $X$, $Y$ are
topological spaces, is said to be faintly precontinuous if for
each    open-closed set $W \subseteq Y$ the set $f\sp {-1}(W)$
is preopen. The paper contains: characterizations of faintly
precontinuous functions, relations to other classes of functions,
theorems on preservation of faint precontinuity under operations
such as superpositions, restrictions, products. Also the
invariance of some separation axioms and properties similar to
compactness are considered. Various examples are given. 

$\bullet$ Jin Han Park and Yong Beom Park, On $sp$-regular
spaces, {\em J. Indian Acad. Math.},  {\bf 17} (2) (1995),
212--218. In this paper, the authors call a topological space $X$
$sp$-regular if, for each closed set $F$ and point $x\in X-F$,
$F$ and $x$ are contained in disjoint semi-preopen sets. They
obtain some equivalent conditions, and study subspaces, products
and continuous images of $sp$-regular spaces. 

$\bullet$ S.~Raychaudhuri and M.N.~Mukherjee, Further
characterizations of $\delta\sb p$-closed spaces, {\em J. Pure
Math.}, {\bf 12} (1995), 27--35) continued of the study of
$\delta\sb p$-closedness of topological spaces, a sort of
covering property strictly stronger than quasi $H$-closedness and
independent of compactness. 

$\bullet$ U.D.~Tapi and S.S.~Thakur and A.~Sonwalkar,
Quasi-preopen sets,
{\em J. Indian Acad. Math.}, {\bf 17} (1) (1995), 8--12. 

$\bullet$ H.~Maki, J.~Umehara and T.~Noiri, Every topological
space is pre-$T\sb {1/2}$, {\em Mem. Fac. Sci. K"chi Univ. Ser.
A Math.}, {\bf 17} (1996), 33--42. A subset $A$ of $(X,\tau)$ is
said to be generalized closed if ${\rm cl}(A) \subseteq U$
whenever $A \subseteq U$ and $U$ is open in $(X,\tau)$. A subset
$A$ of $(X,\tau)$ is said to be pre-generalized closed (denoted
pg-closed) if ${\rm pcl}(A) \subseteq U$ whenever $A\subset
U$ and $U$ is preopen in $(X,\tau)$. A space $(X,\tau)$ is said
to be pre-$T\sb {1/2}$ if a subset is pg-closed if and only if
it is preclosed. In this paper, the authors introduce the concept
of a pre-generalized closed set and prove that every topological
space is pre-$T\sb {1/2}$. This is a strengthening of a result
by A.~Kar and P.~Bhattacharyya (Bull. Calcutta Math. Soc. {\bf
82} (5) (1990), 415--422) that every topological space is
pre-$T\sb 0$. 

$\bullet$ In (D.~Somasundaram and V~Padmavathy, On
generalizations of $H$-closed spaces, {\em Indian J. Pure Appl.
Math.}, {\bf 27} (6) (1996), 557--573.) The authors generalised
$H$-closed spaces by using preopen sets. Analogues of the
$S$-Uryson closed spaces and $s$-regular closed spaces studied
by S.P.~Arya and M.P.~Bhamini ({\em Indian J. Pure Appl. Math.},
{\bf 15} (1) (1984), 89--98) are considered. An adaptation of the
methods employed for HP-closed spaces by T.G.~Raghavan and
I.~Reilly ({\em Indian J. Math.}, {\bf 28} (1) (1986), 75--88)
is used to study projective maximum and projective minimum. 

$\bullet$ D.~Andrijevi\'c, On $b$-open sets, {\em Mat. Vesnik},
{\bf 48} (1-2) (1996), 59--64. The author introduced the class
of $b$-open sets, which is contained in the class of semi-preopen
sets and contains all semi-open sets and all preopen sets. It is
proved that the class of $b$-open sets generates the same
topology as the class of preopen sets.

$\bullet$ E.~Hatir, T.~Noiri and S.~Yksel, A decomposition of
continuity, {\em Acta Math. Hungar.}, {\bf 70} (1-2) (1996),
145--150. The authors define a set $S$ to be a $\cal C$-set if
$S = U \cap A$, where $U$ is open and $A$ is semi-preclosed. The
authors show that a subset of a topological space is open if and
only if it is an $\alpha$-set and a $C$-set. They define a new
class of functions between topological spaces as follows. A
function $f\colon X\to Y$ is defined to be $C$-continuous if the
preimage (under $f$) of each open set in $Y$ is a $C$-set in $X$.
This enables them to provide the following decomposition of
continuity: a function is continuous if and only if it is
$\alpha$-continuous and $C$-continuous. 

$\bullet$ Recently, \'A.~Cs sz r, (Generalized open sets, {\em
Acta Math. Hungar.}, {\bf 75} (1-2) (1997), 65--87) provided a
unified approach to the study of generalized open (in particular
preopen) sets in a topological space, by using monotonic
operators $\gamma \colon \exp(X) \rightarrow \exp(X)$. 

\
E-mail: {\tt dontchev@cc.helsinki.fi}
\
\end{document}